\documentclass[11pt]{article}

\usepackage{graphics}
\usepackage{color}

\usepackage[all]{xy}

\newfont{\frak}{eufm10 scaled\magstep1}
\newfont{\sfrak}{eufm8 scaled\magstep1}
\newfont{\bbb}{msbm10 scaled\magstephalf}
\newfont{\sbbb}{msbm7 scaled\magstephalf}

\def\cdd{\C^d_{\Delta}}
\def\pic{\pi_{\SCC}}
\def\Pic{\Pi_{\SCC}}

\def\cquo{\cdd/\Nc}

\def\dc{\d_{\SCC}}

\def\Nc{N_{\SCC}}

\def\D{\Delta}
\def\C{\mbox{\bbb{C}}}
\def\R{\mbox{\bbb{R}}}
\def\Z{\mbox{\bbb{Z}}}
\def\K{\mbox{\bbb{K}}}
\def\Q{\mbox{\bbb{Q}}}

\def\cd{\C^d}

\def\vz{\underline{z}}

\def\vw{\underline{w}}

\def\ed{e_1,\ldots,e_d}
\def\ld{\lambda_1,\ldots,\lambda_d}
\def\xd{X_1,\ldots,X_d}
\def\d{\mbox{\frak d}}
\def\n{\mbox{\frak n}}
\def\ddu{\d^*}

\def\zd{(z_1,\cdots,z_d)}
\def\lorw{\longrightarrow}
\def\SCC{\mbox{\sbbb{C}}}

\def\s{\mbox{\frak{s}}}

\def\b{\mbox{\frak{b}}}

\def\G{\Gamma}

\def\s{{\mbox{\frak{s}}}}
\def\spic{{\mbox{\sfrak{s}}}}

\def\b{{\mbox{\frak{s}}}}

\def\squareforqed{\hbox{\rlap{$\sqcap$}$\sqcup$}}
\def\qed{\ifmmode\else\unskip\quad\fi\squareforqed}
\def\smartqed{\def\qed{\ifmmode\squareforqed\else{\unskip\nobreak\hfil
\penalty50\hskip1em\null\nobreak\hfil\squareforqed
\parfillskip=0pt\finalhyphendemerits=0\endgraf}\fi}}

\newtheorem{thm}{Theorem}[section]

\newtheorem{lemma}[thm]{Lemma}

\if@titlepage
  \newenvironment{resume}{%
      \titlepage
      \null\vfil
      \@beginparpenalty\@lowpenalty
      \begin{center}%
        \bfseries \resumename
        \@endparpenalty\@M
      \end{center}}%
     {\par\vfil\null\endtitlepage}
\else
  \newenvironment{resume}{%
      \if@twocolumn
        \section*{\resumename}%
      \else
        \small
        \begin{center}%
          {\bfseries \resumename\vspace{-.5em}
\vspace{9pt}}%
        \end{center}%
        \quotation
      \fi}
      {\if@twocolumn\else\endquotation\fi}
\fi

\newcommand\resumename{R\'esum\'e}

\title{\sc Complex quotients by nonclosed groups and their stratifications
}
\author{\sc Fiammetta Battaglia\thanks{Partially supported 
by  GNSAGA (CNR)}
}
\date{}
\begin{document}
\maketitle
\begin{abstract}
We define the notion of complex stratification by quasifolds
and show that such stratified spaces occur as complex quotients by certain
nonclosed subgroups of tori associated to convex polytopes.
The spaces thus obtained provide a natural generalization, to the nonrational
case, of the notion of toric variety associated with a rational convex polytope.
\end{abstract}

\begin{resume}
On d\'efinit la notion de stratification complexe de quasifolds et on
montre que ces espaces stratifi\'es se r\'ealizent comme quotients complexes 
par de sousgroupes non ferm\'es de tori, associ\'es aux polytopes convexes.
Les espaces ainsi obtenus donnent une g\'en\'eralization naturelle, au cas non 
rationnel, de la notion de vari\'et\'e torique associ\'ee \`a un polytope 
convexe rationnel.
\end{resume}

\medskip

{\small 2000 \textit{Mathematics Subject Classification.} Primary: 14M25.
Secondary: 53D20, 32S99, 32C15.}


\section*{Introduction}
Our aim is to define a geometric object that
naturally generalizes, to the nonrational setting, the notion
of toric variety associated with a rational convex polytope.
Let $\d$ be  a real vector space of dimension $n$ and let $\D$ be a
polytope in $\ddu$, rational with respect to a lattice $L$, with
$d$ faces of codimension $1$. Then the toric variety corresponding to $\D$ is
the categorical quotient of
a suitable open subset of $\C^d$, modulo the action of a subtorus of
$(\C^*)^d$, as shown by Cox in \cite{cox}. On the other hand
complex toric spaces corresponding to {\em simple} convex polytopes, not
necessarily rational, were
constructed as complex quotients in joint work with 
Elisa Prato \cite{cx}, who had previously given the
construction of these spaces as symplectic quotients
\cite{p}. More precisely,
to each (simple) convex polytope in $\d^*$, the moment polytope, 
there correspond a unique fan in $\d$--generated by the $1$-dimensional
cones relative to the codimension-$1$ faces of the polytope--and 
many choices of the following data: the generators
of the $1$-dimensional cones of the fan 
and a quasilattice $Q$ in $\d$ containing such generators.
A quasilattice in the vector space $\d$ is a $\Z$-submodule of 
$\d$ generated by a set of generators of $\d$ (cf. for example
\cite{S}). In \cite{cx}, for each choice of data relative to the polytope 
$\D$, an $n$-dimensional toric space is contructed,  
given by the geometric quotient of an
open subset $\C^d_{\D}$ of $\C^d$ by the action of a non necessarily closed
subgroup $N_{\SCC}$ of $(\C^*)^d$. 
The space thus obtained has the
structure of a complex quasifold and is endowed with the holomorphic action
of the $n$-dimensional complex quasitorus $D_{\SCC}=\d_{\SCC}/Q$.
Quasitori are a first, natural example of
quasifolds; they have been introduced, together with the notion of quasifold,
in Prato's article \cite{p}. Complex quasifolds are topological spaces 
whose local models are open subsets of $\C^n$ modulo the action of finitely 
generated groups. Local models are glued together to give rise to a
global quasifold structure. Quasifolds structures are highly
singular: for example the finitely generated groups obtained are in general 
nonclosed, so that the corresponding spaces are non Hausdorff.

Let us now consider a {\em nonsimple} convex polytope, not necessarily
rational;
once generators of $1$-cones and a quasilattice are chosen, we
carry on the construction of our space as complex quotient.
Both $\C^d_{\D}$ and $N_{\SCC}$
can still be defined, the open set $\C^d_{\D}$ depends only on the
combinatorics of the polytope, or, equivalently, of the associated fan, 
whilst the group $N_{\SCC}$ depends on the choice of generators and 
quasilattice.
The first problem we have to deal with is to make sense of
the quotient $\C^d_{\D}//N_{\SCC}$:
as in the rational nonsimple case, we cannot simply take the geometric
quotient, hence 
in Section~2 we define a suitable notion of quotient. Then,
in Section~3, we work out its structure.
Consider a $p$-dimensional face $F$ of $\D$ and a point $\nu\in F$.
The polytope $\D$, in a
neighborhood of $\nu$, is given by the product of $F$ by a cone over an
$n-p-1$ dimensional polytope $\D_F$. We denote by $X_{\D_F}$ the complex space associated to $\D_F$, together with the
induced choice of generators
of $1$-cones and of quasilattice. The face $F$ is regular if
$\D_F$ is a simplex, singular otherwise.
As in the rational case, each
$p$-dimensional face of the polytope corresponds to a $p-$dimensional orbit of
the quasitorus $D_{\SCC}$ acting on $X_{\D}$. In particular the interior
of the polytope corresponds to the dense open orbit of $D_{\SCC}$ in
$X_{\D}$. Orbits produce a stratification of the quotient
$X_{\D}$ that mirrors
the structure of the associated polytope $\D$.
The union of orbits corresponding to regular faces gives the
regular set of $X_{\D}$. Orbits that corresponds to singular faces are the
singular strata. Intuitively it is clear that each stratum has a
natural structure of complex quasifold,
since it is an orbit (or union of orbits) of the quasitorus $D_{\SCC}$.
 Furthermore, in a neighborhood of
a singular orbit, the quotient $X_{\D}$ is the product, in a suitable way,
of the singular stratum by a complex cone over $X_{\D_F}$;
it is a twisted product by a finitely generated group. This group 
is finite in the rational case, when the decomposition
is indeed locally trivial and strata are smooth, namely the usual
notion of stratification is satisfied (cf. \cite{GMcP}).
The definitions of complex cone and complex stratification are given 
in Section~1.
We can then state our main result:
given a convex polytope $\D$,
to each set of generators and of quasilattice containing them,
there corresponds a complex quotient $X_{\D}$, 
endowed with an $n$-dimensional complex stratification by quasifolds and
acted on holomorphically
by a complex quasitorus
with a dense open orbit.
Moreover the space $X_{\D}$ is homeomorphic to
its symplectic counterpart, obtained as symplectic quotient (cf. \cite{bp2,f}).
Such homeomorphism respects
the decompositions and its restriction to each stratum is a
diffeomorphism, with respect to which
the symplectic and complex structures of strata are compatible.
In particular
the space $X_{\D}$ is compact and its strata are in fact Ka\"hler quasifolds.
Remark that, as in the simple case (cf. \cite{p},\cite{cx}), 
the space $X_{\D}$ as complex space depends only on the set of generators and 
on the quasilattice, while its symplectic structure depends on the polytope.

These results complete the generalization of the notion of toric spaces
associated to nonrational convex polytopes, in particular they shed light onto
the relationship with the theory of classical toric varieties.
The geometry and topology of our spaces and
the relationship with the properties of the polytope are  natural
questions related to our work. A first step towards a better understanding
of these different aspects,  that will be pursued in the sequel, is to
define and investigate cohomological invariants  of our spaces.
In this note we give the general idea of our
construction and results, leaving the details of the proofs to the forthcoming
paper \cite{future}.

\vspace{.3cm}

{\bf Acknowledgments}. We would like to thank the referee for useful 
remarks.

\vspace{.3cm}

\section{Complex stratifications by quasifolds}

For the detailed definitions of real symplectic quasifold, quasitorus
and
related notions we refer the reader to \cite{p}, for the complex version
of these
notions see \cite{cx}.
Roughly speaking, as we have already observed, a complex quasifold of
dimension $n$ is locally modelled  on the topological quotient of an open
subset of $\C^n$ by the action of a finitely generated group.
A basic example of real quasifold is Prato's
quasicircle $D^1_{\alpha}=\R/\Z+\alpha\Z$
with $\alpha\in\R\setminus\Q$, this gives also an example of  quasitorus
(cf. \cite{p}). Notice that, if $\alpha$ is taken in $\Q$, then $D^1_{\alpha}$
is either an orbifold or $S^1$.
The complexification
$\left(D^1_{\alpha}\right)_{\SCC}$ of $D^1_{\alpha}$ is the complex
quasitorus given by $\C/\Z+\alpha\Z$.

The notions of decomposition and stratification of a space
$X$ were given in \cite{f}.
We allow the pieces of a decomposition to be quasifolds and
we require then the usual properties of a decomposition.
The dimension of the maximal piece is by definition the dimension of $X$, 
say $n$.
A smooth map (resp. an isomorphism) between decomposed spaces is a
continuos map (resp. homeomorphism) that
respect the decompositions and is smooth (resp. a diffeomorphism)
when restricted to pieces.
A stratification is a decomposition that locally, near to
each point $t$ of an $r$-dimensional stratum $\cal T$,
is given by a twisted product of the kind
$\tilde{B}\times C(L)/\Gamma$, where
$\tilde{B}/\Gamma$ is a local model of
$\cal T$
around
$t$ such that the finitely generated group $\G$ acts freely on
$\tilde{B}$; $L$ is
a $(n-r-1)$-dimensional compact space,
called the {\it link} of $t$, decomposed by
quasifolds and endowed with an action of $\Gamma$ that
preserves the decomposition; $C(L)$ is a real cone over $L$.
Then the decomposition of the link $L$  is required to
satisfy the
above condition. Recursively
we end up, after a finite number of steps,
with links that are compact quasifolds.

We then give the notion
of complex structure of the stratified space.
Our requirements are very strong and are not usually satisfied by complex
stratified spaces, for example the local trivialization is usually
far from being
holomorphic, however they are verified by toric varieties as well as
by our toric spaces. More precisely we require that
for each link $L$ there exists a compact space $Y$ decomposed
by quasifolds and a smooth surjective map $s: L\lorw Y$,
with fibers diffeomorphic  to a fixed $1$-parameter subgroup $S$ of a real 
torus;
each piece of the decomposed spaces $X$, $C(L)$'s and 
$Y$'s is endowed with a complex structure and
the natural projection $C(L)\setminus\{\hbox{cone pt}\}\lorw Y$,
induced by $s$, when restricted to each piece, is holomorphic,
with fibers biholomorphic to the complexification of $S$; moreover
the space $X$ is locally biholomorphic to  the product
$\tilde{B}\times C(L)/\Gamma$, that is the identification mapping
is not only a diffeomorphism but a biholomorphism when restricted to pieces.
We call $C(L)$ a complex cone over $Y$.

\section{Complex quotients by nonclosed groups}
Our goal here is to generalize the construction given by Cox:
we need to produce, in association to
$\D$, an open subset of $\C^d$ and a subgroup of the torus
$(\C^*)^d$.
We first take care of the open subset.
Consider the open faces of $\Delta$. They can be described
as follows.
Write the polytope as intersection of half spaces:
$
\D=\bigcap_{j=1}^d\{\;\mu\in\ddu\;|\;\langle\mu,X_j\rangle\geq\lambda_j\;\}
$
for inward pointing vectors $\xd$ in  $\d$
and the real
numbers $\ld$ determined by our choice of $X_i$. Notice that the vectors $X_i$ are
generators of the $1$-dimensional cones of the fan in $\d$ dual to the
polytope.
For each face $F$ there exists a
subset $I_F\subset\{1,\ldots,d\}$ such that
$F=\{\,\mu\in\Delta\;|\;\langle\mu,X_j\rangle=\lambda_j\;
\hbox{ if and only if}\; j\in I_F\,\}.$
The $n$-dimensional open face of $\D$ corresponds to the empty subset.
Define the open set
$\widehat{V}_F=\{\,\zd\in\cd\;|\;z_j\neq0\;\;\hbox{if}\;\; j\notin I_F\,\}$
and  denote by $\C^d_{\D}$ the open subset of $\C^d$ given by
$\C^d_{\D}=\cup_{F\in\D}\widehat{V}_F.$
Notice that, in the definition of the open subset $\C^d_{\D}$,
only the combinatorics of the polytope intervenes. Moreover, the open subset
$\C^d_{\D}$ coincides with the one defined in \cite{cox} for the rational
case. It is in the definition of the group acting on $\C^d_{\D}$ that
nonrationality intervenes. In order to define the 
group $N_{\SCC}$
we adopt the following procedure: it is an extension of the procedure
introduced by Delzant in \cite{delzant},  extended to the nonrational case
first in \cite{p} and then in \cite{cx,f}.
Let us fix a quasilattice $Q$ in the space $\d$
containing the elements $X_j$ (for example
$\hbox{Span}_{\Z}\{X_1,\cdots,X_d\}$).
Consider the surjective linear mappings $\pi \,\colon\,\R^d \lorw \d$
(respect. $\pic \,\colon\, \C^d  \lorw  \dc$) defined by
$\pi(e_j)=X_j$ (respect. $\pic(e_j)=X_j$),
with $\{\ed\}$ the standard basis.
Consider the quasitorus $\d/Q$ and its complexification
$\dc/Q$. The mappings $\pi$ and $\pic$ induce the group homomorphisms
$\Pi
\,\colon\, (S^1)^d=\R^d/\Z^d\lorw \d/Q$ and $\Pic \,\colon\,
(\C^*)^d=\C^d/\Z^d\lorw \dc/Q$ respectively.
We define $N$ (respect. $\Nc$) to be the kernel of the mapping $\Pi$
(respect. $\Pic$). The group $N$ has dimension $(d-n)$.
Let $\n=\hbox{Ker}\,\pi$ the Lie algebra of $N$, then
for the complexified group $\Nc$
the polar decomposition holds, namely
\begin{equation}\label{polareq}
\Nc=NA,
\end{equation}
where $A=\exp(i\n)$. 
If $Q$ is an honest lattice then $N$ is a compact real torus.

For a given polytope
a choice of normals and of quasilattice $Q$ is said to be {\it rational}
if  $Q$ is a true lattice.
There are polytopes that do not admit rational choices and
that are not combinatorially equivalent to rational polytopes \cite{gru}.
If the polytope is rational in a lattice $L$
the categorical quotient
$\C^d_{\D}//N_{\SCC}$, constructed by Cox in \cite{cox},
can be thought of as the quotient
of $\C^d_{\D}$ by the following equivalence relation:
two points in $\C^d_{\D}$ are
equivalent if the closures of the $N_{\SCC}$-orbits
through these points have nonempty
intersection.
In our context this definition does
not lead anywhere, since the group $N$ itself is nonclosed. Therefore we have
to distinguish two different ways in which orbits are non closed, the first is
given by the fact that $N$ is nonclosed, this is peculiar to the nonrational
setting and produces the quasifold structure of strata, the other is due to the
fact that, as in the rational case, there are {\em nonclosed} $A$-orbits.
There is absolutely no difference
of behavior when it comes to
$A$-orbits, as an example consider $e^{2\pi a t}$,
with $t\in\R$ and $a$
a real constant. Let $z\in\C$, obviously the orbit  $e^{2\pi a t}z$
is not influenced by having $a$ rational or not.
Therefore we are led to consider $A$-orbits.
Let $\vz\in\C^d_{\D}$, we say that the $A$-orbit $A\vz$ is closed if it is
closed in $\C^d_{\D}$. Let $J$ be any set of indices in $\{1,\dots,d\}$, we
denote by $\K^J=\{\vz\in\C^d\;|\;z_j\in\K\quad\hbox{if}\quad j\in J,
z_j=0\quad\hbox{if}\quad j\notin J\}$, where $\K$ can be either $\C$ or
$\C^*$.
Consider the $(\C^*)^d$-orbit
$(\C^*)^{I_F^c}=\{\,\zd\in\cd\;|\;z_j=0\;\;\hbox{iff}\;\; j\in I_F\,\}.$
\begin{thm}\label{closedorbits}
Let $\vz\in\C^d_{\D}$. Then the $A$-orbit through $\vz$, $A\vz$, is
closed if and only if there exists a face $F$ such that $\vz\in
(\C^*)^{I_F^c}$.
Moreover, if $A\vz$ is nonclosed, then it contains one and only one
closed $A$-orbit.
\end{thm}
Theorem~\ref{closedorbits}, which holds in the rational setting too,
allows us to define on the open set $\C^d_{\D}$ the
following equivalence relation: two points $\vz$ and $\vw$ are equivalent
if and only if $$\left(N(\overline{A\vz})\right)\cap\left(\overline{A\vw}
\right)\neq\emptyset$$
where the closure is meant in $\C^d_{\D}$.
We define the space $X_{\D}$ to be the quotient of 
$\C^d_{\D}$ by the equivalence
relation just defined, we denote the quotient by
$X_{\D}=\C^d_{\D}//N_{\SCC}.$
Notice that, if the polytope is simple, then $\C^d_{\D}=
\cup_{F\in\D}(\C^*)^{I_F^c}$
and
the quotient $X_{\D}$ is just the geometric quotient,
whilst if the polytope is nonsimple and rational the quotient 
$X_{\D}$ is the known categorical quotient.

\section{The stratification}

The decomposition in pieces of the quotient $X_{\D}$ reflects the
geometry of the polytope $\D$.
A partial order on the set of all faces of $\D$ is defined by
setting $F\leq F'$ if
$F\subseteq\overline{F'}$. The polytope $\Delta$ is the disjoint
union of its faces. Let $F$ be a $p$-dimensional face and let
$r_F=\hbox{card}(I_F)$, clearly
$r_F\geq n-p$. When the equality holds
the face $F$ is said to be {\em regular}, {\em singular} otherwise.
Consider $\d_F=\hbox{Span}\{X_j\;|\;j\in I_F\}$, the
natural injection $j_F:\d_F\hookrightarrow \d$ and the subset
$\Sigma_F^{\diamond}=\cap_{j\in I_F}\{\xi\in\d^*\;|\;\langle
\xi,X_j\rangle\geq\lambda_j\}$. Then
$j_F^*(\Sigma_F^{\diamond})=\Sigma_F$ is a polyhedral cone with 
vertex $j_F^*(F)$.
By cutting this cone with an affine hyperplane, transversal to its
codimension 1 faces, we obtain an $(n-p-1)$-dimensional polytope,
that we call $\D_F$, which of course depends
on the choice of the hyperplane.
Each $q$-dimensional face $G$ of $\D$
greater than $F$ gives a face in $\D_F$ of dimension $q-p-1$, singular
if and
only if $G$ is singular in $\D$. Near to $F$ the polytope
$\D$ is the product of $F$ by a cone over $\D_F$, and this is exactly the
stratified structure of the toric space that we have constructed.
The maximal stratum is
${\cal T}_{\scriptstyle{max}}=
\cup_{F\;\scriptstyle{reg}}((\C^*)^{I_F^c})/N_{\SCC},$
whilst there is a piece ${\cal T}_F$ for each singular face $F$
of $\D$, which can be identified with
the orbit space
${\cal T}_F=(\C^*)^{I_F^c}/N_{\SCC}.$
Theorem~\ref{closedorbits} implies that
the space $X_{\D}$ is given by the union of the strata just defined.
The structure of $X_{\D}$ as decomposed space and the properties that 
characterize
$X_{\D}$ as toric space associated to $\D$ are given in the following
statement:
\begin{thm}\label{stratiquasifold}
The subset ${\cal T}_F$ of $X_{\D}$ corresponding to each $p$-dimensional
singular face of $\D$ is a $p$-dimensional complex quasifold.
The open subset ${\cal T}_{max}$  
is an $n$-dimensional complex quasifold. These subsets give a
decomposition by complex quasifolds of $X_{\D}$. Moreover there is a continuos
action of $D_{\SCC}$ on $X_{\D}$, with the dense open orbit $(\C^*)^d/N_{\SCC}$. Such action is holomorphic when restricted to pieces.
\end{thm}
The maximal stratum is said to be the regular stratum, 
it is the analogue of the open set of rationally smooth points in the rational 
case, whilst the strata corresponding to singular faces are singular (see
Lemma~\ref{prodotto}).

Now consider the cone $\Sigma_F$ in the space $\d^*_F$, 
with induced normal vectors $X_j\in\d_F$, 
$j\in I_F$, and quasilattice $Q\cap \d_F$. 
Then let $X_0=\sum_{j=1}^d s_jX_j$, with
$s_j\in(0,1)$ 
suitably chosen for $j\in I_F$ and 
$s_j=0$ for $j\notin I_F$. We consider $X_0$ in $\d_F$ and we
denote by
$\hbox{ann}(X_0)$ the annihilator 
of $X_0$. We can view $\D_F$ as a convex
polytope lying in the linear hyperplane $\hbox{ann}(X_0)
\stackrel{k_F}{\hookrightarrow}\d_F^*$,
with induced normal vectors
$k_F^*(X_j)$ in $\hbox{ann}(X_0)^*\simeq \d_F/\langle X_0 \rangle$,  
$j\in I_F$, and quasilattice $k_F^*(Q\cap \d_F)$.
The spaces corresponding to the convex sets $\Sigma_F$ and $\D_F$ are 
$X_{\Sigma_F}=\C^{I_F}//N^F_{\SCC}$
and $X_{\D_F}=\C^{I_F}_{\D_F}//(N^F_0)_{\SCC}$,
where the $(r_F-n+p)$-dimensional group
$N^F_{\SCC}$ is given by $N_{\SCC}\cap (\C^*)^{I_F}$
and
$(N^F_0)_{\SCC}/N^F_{\SCC}\simeq \exp(\s+i\s)$, with
$\s=\hbox{Span}\{(s_1,\dots,s_d)\}$.
We first prove that the stratification satisfies the local triviality
condition:
\begin{lemma}\label{prodotto} Let $F$ be a singular face, then
the singular stratum ${\cal T}_F$ can be identified with $(\C^*)^p/\G_F$,
where $\G_F$ is a finitely generated subgroup of $(\C^*)^p$ acting
on $X_{\Sigma_F}$ and
freely
on $(\C^*)^p$. There is a mapping from  
$(\C^*)^{p}
\times X_{\Sigma_F}/
{\G}_{F}$ onto the open subset $(\C^{I_F}\times(\C^*)^{{I_F}^c})
//N_{\SCC}$ of $X_{\D}$, which is a homeomorphism and a  biholomorphism when
restricted to the pieces of the respective decompositions.
\end{lemma}
Now, in order to complete the proof that the space $X_{\D}$ is a stratified 
space, it remains to show that:
\begin{itemize}
\item[($\ast$)] 
for each singular face $F$ in $\D$  there is a link $L_F$ which satisfies
the definition, moreover the cone $X_{\Sigma_F}$, the link $L_F$ and the toric
space $Y_F=X_{\D_F}$
have the properties required for a stratification to be complex.
\end{itemize}
General results, for example on bundles, are often not readily applicable 
to our spaces because of their topology, therefore,
although a description of the spaces in statement ($\ast$) can be given 
within our set up (see the first row of diagram~\ref{diagramma1}), we make 
use of
the interplay with the symplectic quotients in order to give a neat description
of the link $L_F$ and of $X_{\Sigma_F}$ as a real cone over it. Hence,
before going on to describe the cone $X_{{\Sigma}_F}$, 
let us briefly recall from \cite{f}
the symplectic construction. Let $\Psi_{\D}:\C^d\lorw\n^*$ be the
moment mapping with respect to the $N$-action such that
$\Psi_{\D}(0)=\sum_{j=1}^d\lambda_j\iota^*(e_j^*)$, where
$\iota:\n\lorw \R^d$ is the inclusion.
The quotient
$M_{\D}=\Psi_{\D}^{-1}(0)/N$ is endowed with a symplectic stratification by
quasifolds and is acted on effectively by the quasitorus $D=\d/Q$. Moreover
there is a continuos mapping $\Phi:M_{\D}\lorw\d^*$ whose restriction
to each stratum, with its symplectic structure, is a moment mapping with
respect to the action of $D$, the image of $\Phi$ is exactly $\D$.
The inclusion $\Psi^{-1}_{\D}(0)\hookrightarrow\C^d$ induces a continuos
mapping $\chi_{\D}:\Psi^{-1}_{\D}(0)/N\lorw \C^d_{\D}//N_{\SCC}$.
It was proved in \cite{cx} that, when the polytope is {\em simple}, the mapping
$\chi_{\D}$ is a diffeomorphism inducing a Kahler structure on $X_{\D}$.
For any convex polytope, it remains to prove that:
\begin{itemize}
\item[($\diamond$)] 
the mapping $\chi_{\D}$ identifies the symplectic and complex quotients
as stratified spaces.
\end{itemize}
Statements ($\ast$) and ($\diamond$)  are stricly intertwined and
they can 
be 
proved 
together by induction on the depth of the polytope $\D$, which is the maximum 
length that a chain of singular faces can attain in $\D$.
The key diagram is the following
\begin{equation}\label{diagramma1}
\xymatrix{
{\scriptstyle
\C^{I_F}//(N^{F}_{\SCC})\setminus\{
 [0]
\}}
\ar[r]^{q_2}
\ar@/^2pc/[rr]|s
\ar @{} [dr]|*+[o][F-]{2}
&
{\scriptstyle
\C^{I_F}_{\D_{F}}//N^{F}_{\SCC}\exp(i\spic)}
\ar @{} [dr]|*+[o][F-]{1}
\ar [r]^{q_1}
&
{\scriptstyle \C^{I_F}_{\D_{F}}//{N^{F}_0}_{\SCC}}
\\
{\scriptstyle
(\Psi_{\Sigma_F}^{-1}(0)/N^{F})\setminus \{[0]\}}
 \ar[u]^{\chi^2_{F}}\ar[r]^{p_2}
\ar@/_2pc/[rr]|{s'}
&
\scriptstyle{
(\Psi_{\D_F})^{-1}(0)/N^{F}}
 \ar[u]^{\chi^1_{F}}\ar[r]^{p_1}&
\scriptstyle{
(\Psi_{\D_F})^{-1}(0)/N^{F}_0}
 \ar[u]^{\chi^0_{F}}
}
\end{equation}
The mappings $\chi_F^j$ are all diffeomorphisms of decomposed spaces,
both diagrams commute, the projections $q_2,p_2$
allows us to identify the space $X_{\Sigma_F}$ to a real cone
over the link $L_F$, which is the space
in the mid column. The projections $p_1,q_1$ are fibrations of the link
$L_F$ over the compact Kahler space $X_{\D_F}$, with fibre $\exp(\b)$,
the projections $s,s'$ are fibrations of the complex space
$X_{\Sigma_F}\setminus\{\hbox{cone pt}\}$
over the compact space $X_{\D_F}$, with fibre $\exp(\b+i\b)$. All mappings are
natural, preserve the decompositions and the structure of the strata.
We finally state our main result:

\begin{thm}\label{ultimo}  Let $\d$ be a vector space of dimension $n$, and let
$\D\subset\ddu$ be a convex polytope. Choose inward-pointing normals to the
facets of $\D$, $X_1,\ldots,X_d\in\d$, and let $Q$ be a quasilattice containing
them.  The corresponding quotient $X_{\D}=\C^d_{\D}//N_{\SCC}$ is a 
complex stratified space.
The mapping $\chi_{\D}\,\colon\,\Psi_{\D}^{-1}(0)/N\lorw\cquo$ is an 
equivariant homeomorphism whose restriction
to each stratum is a diffeomorphism of quasifolds, with respect to which the
symplectic  and  complex structure are compatible,
so that strata have the structure
of K\"ahler quasifolds.
\end{thm}

\noindent \small{\sc Dipartimento di Matematica Applicata,
Via S. Marta 3, 50139 Firenze, ITALY,  {\tt mailto:fiammetta.battaglia@unifi.it}}
\end{document}